\newcommand\by[1]{#1,}

\newcommand\paper[1]{\emph{#1},}
\newcommand\jour[1]{#1,}
\newcommand\vvvol[1]{\textbf{#1}}
\newcommand\yr[1]{(#1)}

\newcommand\pages[1]{#1}
\newcommand\book[1]{\emph{#1},}
\newcommand\bookinfo[1]{#1,}

\newcommand\publ[1]{#1,}
\newcommand\publaddr[1]{#1,}

 \documentclass[%
    draft,
     12pt,
   ]%
  {amsart}
\usepackage{%
amssymb, %
latexsym, %
amscd,
amsmath, %
verbatim 
}
\newtheorem{theorem}{Theorem}[subsection]
\newtheorem{corollary}[theorem]{Corollary}
\newtheorem{lemma}[theorem]{Lemma}
\newtheorem{proposition}[theorem]{Proposition}
\theoremstyle{remark}
\newtheorem{remark}[theorem]{Remark}
\newtheorem{remarks}[theorem]{Remarks}

\theoremstyle{remark}

 \newtheorem{notation}[theorem]{Notation}

\numberwithin{equation}{subsection}
\newcommand{\thmref}[1]{Theorem~\ref{#1}}

\newcommand{\lemref}[1]{Lemma~\ref{#1}}
\newcommand{\corref}[1]{Corollary~\ref{#1}}
\newcommand{\propref}[1]{Proposition~\ref{#1}}

\newcommand\begeqn{\begin{equation}}
\newcommand\begeqns{\begin{equation*}}
\newcommand\ndeqn{\end{equation}}
\newcommand\ndeqns{\end{equation*} }
\newcommand\C{\mathbb C}%

\newcommand\detail[1]{}

\newcommand{\eat}[1]{}

\newcommand\Q{\mathbb Q} 

\newcommand\qiii{(qua\-si\kern-.3em~-~\kern-.35em)in\-de\-pend\-ent}
\newcommand\qiind{(qua\-si\kern.25em)\kern-.65em-\kern.35em\-independence}
\newcommand\qiir{(qua\-si\kern-.3em~-~\kern-.35em)re\-la\-tion}

\newcommand\R{\mathbb R} 

\newcommand\T{\mathbb T}%

 \newcommand{\therosteritem}[1]{(#1)}

\newcommand\Z{\mathbb Z} 
\newcommand\zp[1]{Z_{p_{#1}}}

\begin{document}
\title[\texttt{(Draft--\today)} Qua\-si-independence and permutations]
{Permutation and extension for planar
 quasi-independent subsets of the roots of unity}

    \author[Ramsey \& Graham]
{L.~Thomas Ramsey and Colin C. Graham}
\address{Department of Mathematics,
University of Hawaii, Keller Hall, 2565 The Mall,
Honolulu, Ha\-wa\-ii,
96822\\
email\textup{:\texttt{ramsey\@math.hawaii.edu}}}
\address{University of British Columbia.
Mailing address: RR\#1--H-46,
Bow\-en Is\-land, BC, V0N 1G0 Canada.
email\textup{:\texttt{ccgraham@alum.mit.edu}}
}
\thanks{The second named author
is partially supported  by
an NSERC grant.}

\dedicatory{}

\date{\today}

%

\begin{abstract}  Let
$e^{2\pi i\Q}$
denote the set of  roots of
unity.  We consider
subsets $E\subset e^{2\pi i\Q}$
that are quasi-independent or
algebraically independent
(as subsets of the discrete plane).
 A bijective map on $e^{2\pi i\Q}$
  preserves  the
 algebraically independent sets iff
 it preserves the quasi-independent sets,
and those maps are
 characterized.

 The effect on
the size of \qiii\ sets
in the $n^{th}$ roots of unity
$Z_n$
of
increasing a prime factor of $n$
is studied.
\end{abstract}

  \maketitle

             \centerline{\today}

\section{Introduction and statement
of results}\label{SecIntro}

 \subsection{Background}\label{why}

The motivation of this  paper is the
question: {\it what are the properties of those
subsets of $e^{2\pi i\Q}$, the
set of all roots of unity,
that
are quasi-independent?}
That is of interest  for itself
and because a  theorem
of Pisier \cite{pisier} shows  a set $E$ is
a Sidon set iff there is a $\delta>0$ such that
each finite $F\subset E$ contains a quasi-independent
set with $\#F\ge\delta\# E$. Hence, determining the
properties of the quasi-independent subsets of
subgroups of the roots of unity
may help to resolve whether
$e^{2\pi i\Q}$ is Sidon.

The present paper is a continuation
of  \cite{RGS}, where
we studied the size $\Psi(n)$ of the largest
quasi-independent subset of the $n^{th}$
roots of unity and established some
properties of $\Psi$ (some of those results
are quoted in
\S\ref{secPrelims}), and showed (among other things)
that if
a set of primes $P$ has $\sum_{q\in P}1/q<\infty$,
then the set $W$ of roots of unity generated by roots
of the form $e^{2\pi i/(q_1^{n_1}\cdots q_k^{n_k})}$,
where $q_j\in P$ and $1\le n_j$, $1\le j\le k$,
is a Sidon set in $\C$ when $\C$ is given the
discrete topology.

We now define
our terms.

A subset
$B$ of an (additively-written) abelian group
$G$ is
{\it quasi-independent\/} if
\begin{align*}
&k\ge1,\,x_j\in B,\, \epsilon_j=0,\pm1 \text{ for }1\le j\le k \text{ and }\\
&\quad \sum_j\epsilon_j x_j=0\Rightarrow \ \epsilon_j=0
\text{ for } 1\leq j\leq k.
\end{align*}
The set $B$ is {\it independent\/} if
\begin{align*}
&k\ge1,\,x_j\in B,\, \epsilon_j\in\Z\text{ for }1\le j\le k \text{ and }\\
&\quad \sum_j\epsilon_jx_j=0\Rightarrow \epsilon_jx_j=0
\text{ for } 1\leq j\leq k.
\end{align*}
For independence
in divisible groups such as $\C$, it is
equivalent to substitute
``$\epsilon_j \in \Q$''
for ``$\epsilon_j \in \Z$''.
The set $E\subset \C\approx \R^2$ is {\it a Sidon set\/}
if every bounded function on $E$
is the restriction of Fourier-Stieltjes transform of a bounded regular Borel
measure on the Bohr compactification of $\C$.

Throughout this paper $n$ will
denote a positive integer and $T_n$ the
$n$-th roots of unity.
When we use the terms ``quasi-independent'' or ``independent,''
we shall mean as subsets of the additive group $\C $ of complex numbers.
We are interested in interplay between
the additive group
$\C$ and the multiplicative group $\T$,
the set of complex numbers of modulus one,
which is
a closed sub{\it set\/}
of $\C$. We shall ignore the
topology of $\C$. For consistency with
the literature on finite abelian groups
we shall often use additive notation for
the finite subgroups of $\T$ which we study here!
Thus, $T_n$
with multiplication
is identified with the cyclic group $Z_n$
of order $n$ with addition mod $n$.
We often  abuse notation to                                                                                                                                                                                                                                                                                                                                                                                                                                                                                                                                                                                                                                                                                                                                                                                                                                                                                                                                                                                                                                                                                                                                                                                                                                                                                                                                                                                                                                                                                                                                                                                                                                                                                                                                                                                                                                                                                                                                                                                                                                                                                                                                                                                                                                                                                                                                                                                                                                                                                                                                                                                                                                                                                                                                                                                                                                                                                                                                                                                                                                                                                                                                                                                                                                                                                                                                                                                                                                                                                                                                                                                                                                                                                                                                                                                                                                                                                                                                                                                                                                                                                                                                                                                                                                                                                                                                                                                                                                                                                                                                                                                                                                                                                                                                                                                                                                                                                                                                                                                                                                                                                                                                                                                                                                                                                                                                                                                                                                                                                                                                                                                                                                                                                                                                                                                                                                                                                                                                                                                                                                                                                                                                                                                                                                                                                                                                                                                                                                                                                                                                                                                                                                                                                                                                                                                                                                                                                                                                                                                                                                                                                                                                                                                                                                                                                                                                                                                                                                                                                                                                                                                                                                                                                                                                                                                                                                                                                                                                                                                                                                                                                                                                                                                                                                                                                                                                                                                                                                                                                                                                                                                                                                                                                                                                                                                                                                                                                                                                                                                                                                                                                                                                                                                                                                                                                                                                                                                                                                                                                                                                                                                                                                                                                                                                                                                                                                                                                                                                                                                                                                                                                                                                                                                                                                                                                                                                                                                                                                                                                                                                                      abuse notation
identify $T_n$ and $Z_n$ with the
product of cyclic groups of prime  order
and do not distinguish between a factor of those
groups and the isomorphic subgroup.

 \subsection{Statement of  results}\label{results}
 We characterize the permutations
of $Z_n$ that preserve the quasi-independent and
independent subsets of $e^{2\pi i\Q}$,
and show that those classes of
permutations are the same.
The characterization
is of interest in itself, as well as being
helpful in
machine assisted calculations.
 $E \subset \C$
is \qiii\  if and only if the
product  $\rho E$ of $E$ with
a nonzero complex number $\rho$ is \qiii.
Hence quasi-in\-de\-pen\-dence
is invariant
for $E \subset T_n$ under any rotation
permutation of $T_n$, that is,
\qiind\ is invariant
under the group operation in $Z_n$.

From general algebraic principles, we would expect
 quasi-in\-de\-pen\-dence to be invariant
 under group automorphisms.  However,
 there is a much larger set of permutations of
$T_n$   which preserve quasi-indepen\-dence
and independence. Here is a
simplified version of our permutation
result, with the full version being
\thmref{rigidpermuts}.

\begin{theorem}\label{simplerigidpermuts}
Let $n$ be odd and square free, with prime factorization
$n=\prod_{j=1}^K p_j$.
\begin{enumerate}
\item
Let $\sigma_j$ be a permutation of $Z_{p_j}$, $1\leq j\leq K$
and let $\sigma=\sigma_1\times \cdots\times \sigma_K$
be the product permutation
on $Z_n$. Then $\sigma$ preserves both the class of
quasi-independent and the class of
independent sets.
\item
Let $\sigma$ be a permutation of  $Z_n$
that preserves either the class of quasi-independent or the
class of independent sets.
Then $\sigma$ is a product of permutations as in \therosteritem1.
\end{enumerate}
\end{theorem}

The permutations of $e^{2\pi i\Q}$
which preserve quasi-independent
sets are identified in
\corref{permutforexpQ}.

 Our second main result  shows
 that $\Psi$ is ``monotone'' in a certain sense:

\begin{theorem}\label{monotonicity}
For primes $p_1 <p_{2} < \hdots <p_{K}$, let $n=\prod_{j=1}^K p_j$.
Let $s \in \{1, \dots, K\}$.  If $s=K$, let $q$ be any prime such that
$q>p_K$.  If $s<K$, let $q$ be any prime such that $p_s<q$ and $q \notin
\{p_{s+1}, \hdots , p_K\}$.  Let $m=n/p_s$, and
$\Delta\ge0$.
Then the following hold.
\begin{enumerate}
\item
$\Psi(qm)\ge \Psi(n)+(q-p_s)\Psi(m).$
\item
If
$\Psi(n)\ge \phi(n)+\Delta$,
then
$
\Psi(qm) \ge \phi(qm)+\Delta.
$
\item
If
$
\Psi(n) \ge (p_s-1) \Psi(m) +\Delta,
$
then
$
\Psi(qm) \ge (q-1) \Psi(m)+\Delta.
$
\end{enumerate}
\end{theorem}

The factors $p_s-1$ and $
q-1$   are quite natural in this context,
as $\Psi$ has properties quite similar to those of
Euler's $\phi$
\cite[Thm.~4.1]{RGS}.

We give an application of \thmref{monotonicity}
in \S\ref{beyondOneOhFive}.

 \subsection{Organization
 of this paper}\label{org}
Notation and the material needed from \cite{RGS}
 are given in
\S\ref{secPrelims}.
We state and prove the general
permutation theorem in \S\ref{SecPermutation}.
In \S\ref{monotone} we show how $\Psi$
is monotonic and consider extensions of
quasi-independent sets.

When our results are valid for
both the class
of independent sets and for the
class of quasi-independent sets, we write
(quasi-)in\-de\-pend\-ent, (quasi-\nobreak)in\-de\-pend\-ence,
and (qua\-si)\-rel\-a\-tion (see below).

\section{Preliminaries}\label{secPrelims}
\subsection{Simple characterizations of \qiii \ sets}
Giv\-en $n$, it will be convenient to have
specified its prime factorization:
\begin{equation}\label{primefac}
n=\prod_{j=1}^K p_j^{n_j},
\end{equation}
for some $K$ and distinct,
positive primes $p_j$,
and the square-free versions:
\begin{equation}\label{tilden}
\tilde n= \prod_{j=1}^K p_j \text{ and }
Z_{\tilde n}=(\frac n{\tilde n}) Z_n.
\end{equation}

\begin{theorem}\label{primsubgroupscosets}
\cite[Theorem 1.2.2]{RGS}
A set $E\subset  T_n$
 is \qiii \ if and  only if
the intersection of $E$ with each  coset
of $T_{\tilde n} $ is \qiii.
\end{theorem}

A set $E\subset e^{2\pi i\Q}$
is  independent (resp. quasi-independent)
if $E$
supports no non-zero
\emph{relation}
(resp. \emph{quasi-relation}),
that is a  finitely supported
function $f:E\to \Q$ (resp. $E\to \{0,\pm1\}$)
such that
 $\sum_{z\in E}f(z)z=0$.
Switching to $Z_n$, we can describe the
(quasi-)relations via a basis, since the
set of relations is a vector space over $\Q$.

  \begin{lemma}\label{basiszpqr}
  \cite[Lemma 2.10]{RGQ}
 Let $2\le n$ be square-free
 with prime factorization $n=p_1\cdots p_K$.
 Fix $1\le j\le K$ and
 $0\le \ell<p_j$, and write $Z_n=Z_{p_j}\times
 H$. Then  the set of relations on $Z_n$
  has a basis consisting of
  \begin{enumerate}
  \item the characteristic functions
  of cosets of $Z_{p_j}$
  \item relations each of which is
  supported on only one  coset   $H+k$,
  $0\le k<p_j,$ and $k\ne \ell$.
  \end{enumerate}
 \end{lemma}

 Thus, each relation is a sum of
 \therosteritem1 \ ``spikes'' (from cosets of
 $Z_{p_j}$) and \therosteritem2\ of relations supported on
 cosets of $H$, with one coset omitted.

 \begin{corollary}[Empty Floor]\label{smallqiset}
 \cite[Cor.~2.11]{RGQ}
 Let $n, j, H$ be as above.
  Suppose that $E\subset Z_n$ is such that $E\cap (H+b)=\emptyset$
 for some $b$.
 Then
 $E$ is \qiii \ if and only if $E\cap (H+a)$ is \qiii
 \ for all $a\in Z_{p_j}$.
 \end{corollary}

\begin{corollary}\label{psinqge}
Let $n>1$ be square-free and $q$ a positive prime
that does not divide $n$.
Then $\Psi(nq)\ge (q-1)\Psi(n).$
\end{corollary}
\begin{proof}
Let $E\subset Z_n$ be quasi-independent
such that $\#E=\Psi(n)$. Let $F=E\times\{1,\dots,q-1\}
\subset Z_{nq}$. Then $F$ is quasi-independent by
\corref{smallqiset}, and $\#F=(q-1)\#E$.
\end{proof}

\begin{remark}\label{wavyspike}
 Spikes  give us simple examples
 of in\-de\-pen\-dent sets:
let $p$ be a prime factor of the square-free,
odd integer $n$ and $m=n/p$. Suppose that
$E\subset Z_n=Z_m\times Z_{p}$ and $E$
meets each coset of $Z_m$
in at most one point.
Then $E$ is in\-de\-pen\-dent if either $\#E<p$ or
$\#E=p$ but $E$ is not a coset of
$Z_p$.
\end{remark}

 We will need item \therosteritem1 of the following
 theorem.
\begin{theorem}\label{Psiofn}
\cite[Thm.~4.1]{RGS}
For any
integer $n\ge2$, the following hold.
\begin{enumerate}
\item
If a  prime $p$ divides $n$, then
$\Psi(pn)=p\Psi(n)$.
\item If $p$ is prime and $k\ge1$,
then $\Psi(p^k)=\phi(p^k)=p^{k-1}(p-1)$.
\item
If  $n$ is odd, then $\Psi(2n)=\Psi(n)$.
\item Let $p $ be a  prime that does
not divide $n$.
Then $\Psi(pn)\ge (p-1)\Psi(n)$.
\end{enumerate}
\end{theorem}

\begin{theorem}\label{zpq}
\cite[Thm.~6.4]{RGS}
Let $p\neq q$ be odd primes. Then $E\subset Z_{pq}$
is qua\-si-in\-de\-pen\-dent
if and only if all of the following hold. For $t\in Z_p$ and $v\in Z_q$, let
$E^{(t)}= E\cap (Z_q+t)$ and $E_{(v)}=E\cap (Z_p+v)$.
\begin{enumerate}
\item
 The intersection of $E$ with each coset of $Z_p$ is qua\-si-in\-de\-pen\-dent.
 \item
 The intersection of $E$ with each coset of $Z_q$ is qua\-si-in\-de\-pen\-dent.
 \item
For every   $t\in Z_p$ and non-empty subset $F\subseteq E^{(t)}$,
the spikes rising from $F$ are not shadowed by $E$.
\item
For every   $v\in Z_q$ and non-empty subset $F\subseteq E_{(v)}$,
 the spikes rising from $F$
are not shadowed by $E$.
\end{enumerate}
Furthermore, $E$ is in\-de\-pen\-dent if and only if $E$ is qua\-si-in\-de\-pen\-dent.
\end{theorem}

\begin{proposition}\label{minsizeofqiset}
Suppose that $n\ge2$ is  any integer.
Let $p_1$ be the smallest  prime factor of $n$.
\begin{enumerate}
\item
If  $E\subset Z_n$ and $\#E<p_1$, then
$E$ is in\-de\-pen\-dent.
 \item 
 If $E\subset Z_n$ has cardinality $p_1$,
 then $E$ is \qiii \ %
 if and only if $E$ is a not coset of $Z_{p_1}$.
\item 
 If $E\subset Z_n$, $n$ is odd, with at least
 two prime factors, and square-free,
 and $\#E<p_1+p_2-2$,
 then $E$ is \qiii \ if and only if $E$ does
 not contain a coset.
 \item 
 $(Z_{p_1}\cup Z_{p_2})\backslash\{0\}$ is not
 qua\-si-in\-de\-pen\-dent, and has cardinality $p_1+p_2-2$.
\end{enumerate}
\end{proposition}

\begin{proof}
\therosteritem1.
This is trivial when $\#E\le1$ because a non-zero number
in $\C$ is in\-de\-pen\-dent over $\Q$. So we may assume that
$\#E\ge2$.

We argue by induction on $K$, the number of prime factors of $n$.
Suppose  $K=1$. Then for each coset $H$ of
$Z_{\tilde n}=Z_{p_1}$, $\#(H\cap E)<p_1$ so $H\cap E$ is
in\-de\-pen\-dent. By
\thmref{primsubgroupscosets}, $E$ is \qiii.

So suppose $K\ge2$, and that the conclusion
holds for any  $Z_m$ when $m$ has
less than $K$ prime factors.
Suppose that
$n$ has $K$ prime factors $p_1<\dots <p_K$
and that $E\subset Z_n$ has $\#E<p_1$.
Again,  $E$ is
in\-de\-pen\-dent if and only
if $E \cap U$ is in\-de\-pen\-dent
for each coset $U$ of $Z_{\tilde n}$ where
$\tilde n = p_1 \cdot p_2 \cdots p_K.$
Consider any $E \cap U$.  By translation we may
assume $E \subset Z_{\tilde n}$.
We write $Z_{\tilde n}=Z_{p_1} \times H.$

For $0 \le t <p_1$, $E \cap (\{t\} \times H)$
is equivalent by translation to some
 $S \subset H$.  Since $\#(S) \le \#(E)<p_1$
 and $p_1$ is the smallest prime dividing
 the order of $H$, the induction hypothesis
 implies that $S$ and hence $E \cap (\{t\} \times H)$
 are \qiii.
 Since there are $p_1$ disjoint cosets of $H$
 in $Z_{\tilde n}$ but $\#(E)<p_1$, there is some coset of
 $H$ having empty intersection with $E$.
 By \corref{smallqiset}, $E$ is in\-de\-pen\-dent.
 That proves \therosteritem1.

 \smallskip
 \therosteritem2. %
 The ``only if'' direction is clear.
For the converse, suppose that $E$ is not a coset
of $Z_{p_1}$.
If $E$ has elements in different cosets
of $Z_{\tilde n}$, its intersection
with each of those cosets would
have fewer than $p_1$ elements.
By \therosteritem1\ of this Proposition,
each intersection $E \cap U$
would be independent for each
coset $U$ of $Z_{\tilde n}$.
By \thmref{primsubgroupscosets}
\ that would make $E$ be independent.
So we may assume that $E \subset U$
for some coset $U$ of $Z_{\tilde n}$.
By translating we may assume
$U=Z_{\tilde n}$.
We again induct on the number $K$ of prime
factors of $n$.
If $K=1$ and $\#(E \cap U)=p_1$, $E$ would have to be
a coset of $Z_{p_1}$.
So we have $K >1$.

Let $H=Z_{p_1} \times ... \times Z_{p_{K-1}}.$
Since $p_K>p_1$, there are
more than $p_1$ cosets of $H$ within $Z_{\tilde n}.$
Since $\#E=p_1$, at least one coset of $H$ has empty
intersection with $E$.

Suppose $K=2$.  Then $H = Z_{p_1}$ and the intersection
of $E$ within each
coset of $H$ is a proper subset of that coset
(otherwise, $E$ having $p_1$
elements would make it equal that coset).
By \therosteritem1, translated to each coset
of $Z_{p_1}$, the intersection of $E$ with each coset of $Z_{p_1}$ is
in\-de\-pen\-dent.  By \corref{smallqiset}, $E$ is in\-de\-pen\-dent.

 Suppose that \therosteritem2\ holds
for some $K-1 \ge
2$.  Let $H$ be as in the previous paragraph.   Again, since $p_K >p_1$,
there is some coset of $H$ with empty intersection with $E$.
For each coset $W$ of $H$, $E \cap W$ has size at most $p_1$ and $E\cap W$
is not a coset of $Z_{p_1}$.  If  $\#(E \cap W) <p_1$, by the translation of
\therosteritem1\ to the coset $W$, $E \cap W$ is in\-de\-pen\-dent.
If $\#(E \cap W)=p_1$, the induction hypothesis implies that $E \cap W$ is
in\-de\-pen\-dent (if it were a coset of $Z_{p_1}$, then $E=E \cap W$ would be one
also).   By  \corref{smallqiset}, $E$ is in\-de\-pen\-dent.

\smallskip

\therosteritem3. Suppose that $\#E< p_1+p_2-2$
and that $E$ does not contain any cosets.

(i).
Let
$
H=Z_{p_1}\times\cdots
\times Z_{p_{K-1}}\times\{0\}.
$
Suppose $E$ has non-void intersection with every coset of $H$.
By counting,
we see that there is a coset that
meets $E$ in at most one point.\footnote{
In fact,
there at least five such cosets.
 Indeed, let $b$ be the number of cosets with
 exactly one element.
Then  $b + 2(p_K-b) \le p_1 +p_2 - 3$.  Thus
$2 p_K - p_1 -p_2 + 3 \le b$. But  $ p_K \ge p_2$.
Thus $p_K -p_1 + 3 \le b$.   Since
$n$ is odd, $p_K \ge p_2 \ge p_1+2$,
so $5\le b$.}   By translation
(which preserves \qiii \ sets), we may assume that
the zero coset has one point intersection with $E$,
and, by translation again, we may assume that the intersection of
$E$ with $H$ is $\{0\}$.

(ii). Consider any relation supported on $E$,
consisting of the combination of
one spike $f_0$ rising from
$\{0\}$ and relations $f_j$ on the non-zero
cosets $H+j, 1\leq j<p_K$.
The support of $f_j$ is at most $[E\cap (H+j)]\cup [Z_{p_K}\cap (H+j)]$.
 Therefore,
each $f_j=0$, unless $E\cap (H+j)$ has cardinality at
least $p_1-1$. (Throwing in $Z_{p_K}\cap (H+j)$ will
give the remaining point.) If $f_j=0$
for all $j>0$, then $Z_{p_K}\subset E$, contradicting
our assumption that $E$ contains no cosets.

Hence, some $f_j\neq 0$.
We now count. $E$ meets  every
coset of $H$: that's $p_K$ elements. Since some  $f_j\neq0$,
that's at least another $p_1-2$ elements. Hence,
$\#E\ge p_1+p_K-2> p_1+p_2-3$.
Hence, all $f_j=0$. Another contradiction.
That proves that
$E$ cannot meet all cosets of $H$.
We will repeat this argument as follows.

(iii). Let $H+j$ be a coset of $H$ that misses $E$.
By translation, we may assume that $E\cap H=\emptyset$.
 Therefore, $E$ is \qiii \ if and only if
the intersection of $E$ with each of the non-zero cosets of $H$
is \qiii,
by \ref{smallqiset}. Consider  a non-zero coset of $H$.
By translation, we may assume that
the coset is the zero coset. Hence, we have eliminated a factor of $Z_n$,
and $E\subset Z_{p_1}\times \cdots\times Z_{p_{K-1}}$.

(iv). Repeating   paragraphs (i)-(iii), we see that
 we may assume that $K=2$, that is
$Z_n=Z_{p_1}\times Z_{p_2}$.
If $E$ meets all cosets of $H=Z_{p_1}$, then by
paragraphs (i)-(ii),  $E$ contains a coset.
If $E$ misses one coset of $Z_{p_1}$, then
we may assume that $E$ is contained in some other
coset of $Z_{p_1}$, as above, and thus
$E$ is either in\-de\-pen\-dent, or $E$
is (contains) that coset of $Z_{p_1}$.

That proves \therosteritem3.

\smallskip

\therosteritem4. Let the qua\-sirelation $f$ be
definied by $f=\chi_{Z_{p_1}}-\chi_{Z_{p_2}}$.
Then the support of $f$ is
$(Z_{p_1}\cup Z_{p_2})\backslash\{0\}$.
 That proves \propref{minsizeofqiset}.
\end{proof}

\section{The permutation
theorem}\label{SecPermutation}

 \subsection{Preliminary
 remarks}\label{permutssummary}
\begin{notation}\label{preremrigid}
 (i). Let $n\ge2$ be an integer, with
 factorization \eqref{primefac}.
Let $E$ be a full set of representations
of the cosets in $Z_n/Z_{\tilde n}$.
Thus,  $Z_n=Z_{\tilde n}+E$.
\newline
(ii). By ``$Z_{p_a}$'' we remind the reader that
we
mean the subgroup of $Z_n$ of order $p_a$,
so $ 
Z_{p_a}\cong p_a^{n_a-1}Z_{p_a^{n_a}}.
$ 
\end{notation}

\begin{theorem}\label{rigidpermuts}
 Under  \S\ref{preremrigid},  the following hold.
\begin{enumerate}
\item 
Let $1\leq a\leq K$ for some $a$.
Let $\sigma' $ be any permutation of $Z_{p_a^{n_a}}$
that maps each coset of $Z_{p_a}$ into another (or the same)
coset of $Z_{ p_a}$.
Let $\sigma$ be the permutation of $Z_n$ that is the product of
$\sigma'$ with the identity permutations on the other $K-1$ factors
of $Z_n$.
Then $\sigma$ preserves \qiii\ sets.
\item 
Suppose $p_1=2$. Let $\sigma$ be any permutation of $Z_n$
that maps each coset of $Z_{2}$ into itself.
Then $\sigma$ preserves \qiii\ sets.
\item 
Let $\sigma'$ be any permutation of $E$.
Define a permutation $\sigma$ on $Z_n$ by
$\sigma(x+y)=x+\sigma'(y)$ $x\in Z_{\tilde n}, y\in E$.
Then $\sigma$ preserves \qiii\ sets.
\item
Let $\sigma$ be any permutation of $Z_{\tilde n}$.
 Extend
$\sigma$ to all of $Z_n$ by letting $\sigma$
be
the identity on  $Z_n\backslash Z_{\tilde n}$.
Then $\sigma$ preserves the \qiii\ sets of $Z_{n}$
iff it preserves the \qiii \ sets of $Z_{\tilde n}$.
\item  
Any product of permutations of $Z_n$ that preserve the \qiii \ sets
preserves the \qiii \ sets.
\item 
Let $\sigma$ be any group automorphism of $Z_n$ or
a translation of $Z_n$.
Then $\sigma$ preserves \qiii\ sets.
\item 
Every permutation $\sigma$ of $Z_n$ that preserves \qiii \ sets is a product
of permutations of the types of \therosteritem1-\therosteritem4.
\end{enumerate}
\end{theorem}

\begin{remarks}
(i.) The statement of the theorem is slightly redundant,
since \therosteritem3 implies part
of \therosteritem1. Part \therosteritem1 \ is
stated the way it is to assist the reader in understanding the more complicated
\therosteritem3.

(ii). It is easy to see that
\thmref{simplerigidpermuts}
follows immediately from
\thmref{rigidpermuts}.
\end{remarks}

For each prime $p$, let $Z_{p^\infty}$ denote the injective limit of $Z_{p^\ell}$, and
for a finite set of primes $p_1<\cdots<p_K$, let $Z_{n,\infty}$ denote the
product $\prod_1^K Z_{p_j^\infty}$, where
$n=\prod_1^K p_j$.
\begin{corollary}\label{permutforexpQ} Let $\sigma$ be a permutation of
$e^{2\pi i \Q}$
with $\sigma(1)=1$. For each $n\ge2$, we identify $T_n$ with $Z_n$.
Then $\sigma$ preserves the quasi-independent sets iff
$\sigma$ preserves the
independent sets iff for every $n\ge2$ the following hold.
\begin{enumerate}
\item  $\sigma$ maps
each coset of $Z_{\tilde n}$ in $Z_{\tilde n,\infty}$  to another
(or the same)coset of $Z_{\tilde n}$. (The image coset does not have to
be in $Z_{\tilde n,\infty}$.)
\item
If $\sigma(Z_{\tilde n}+t)= Z_{\tilde n}+u$, then there exists a product $\sigma_{\tilde n,t,u}$
of permutations as
in  \thmref{simplerigidpermuts}\ such that $\sigma(x)= \sigma_{\tilde n,t,u}(x-t)+u$
for all $x\in Z_{\tilde n}+t$.
 \end{enumerate}
\end{corollary}
\begin{proof} Straightforward.
 \end{proof}

We will use the following lemma.

\begin{lemma}\label{permutLemma}
Let $\sigma$ be a permutation of $Z_n$ such that,
for $1 \le j \le K$ and $x \in Z_n$, $\sigma(x+Z_{p_j})$ is a coset of $Z_{p_j}$.
Then $\sigma$ preserves \qiind.
\end{lemma}
\begin{proof}  Because $\sigma$ is permutation of $Z_n$, it
sends disjoint cosets of $Z_{p_j}$ into
disjoint cosets of $Z_{p_j}$.  Thus, $\sigma$
induces a one-to-one mapping $\tau$ from $Z_n/Z_{p_j}$
into itself.  Since $Z_n/Z_{p_j}$ is finite,
$\tau$ is a permutation of $Z_n/Z_{p_j}$.

It follows that, for $1 \le j \le K$ and
$x \in Z_n$, $\sigma^{-1}(x+Z_{p_j})$ is
a coset of $Z_{p_j}$.   To see this, observe that
$ 
x+Z_{p_j}=\tau(v+Z_{p_j})=\sigma(v+Z_{p_j}),
$  
\text{for some $v \in Z_n$.}
Since $\sigma$ is a permutation of $Z_n$,
$ 
\sigma^{-1}(z+Z_{p_j})=\sigma^{-1}[\sigma(v+Z_{p_j})]=v+Z_{p_j}.
$ 
Thus, the hypotheses of this Lemma apply
equally to $\sigma$ and $\sigma^{-1}$.

Suppose that $E$ is \qiii.
Consider any (quasi-)relation $f$ supported on
$\sigma(E)$.  Then, with a redundant (non-unique)
choice of coefficients $c_{j,x}$ we have
 $
f=\sum_{j=1}^K \sum_{x \in Z_n} c_{j,x} \chi_{x+Z_{p_j}}.
$ 
 Then
 \begin{equation*}
f\circ \sigma =
  \sum_{j=1}^K \sum_{x \in Z_n}
  c_{j,x} \chi_{\sigma^{-1}(x+Z_{p_j})}.
\]
Note that
\begin{enumerate}
\item
The mapping $g\mapsto g \circ \sigma$ is a linear
algebra isomorphism on $\Q(Z_n)$.
\item
$f\circ\sigma$ is supported on $E$.  To see this, let
$w \in Z_n\backslash E$.  Then, since $\sigma$ is a
permutation of $Z_n$, $\sigma(w) \notin \sigma(E)$
and thus $f(\sigma(w))=0$.
\item
Similarly,
$\chi_{x+Z_{p_j}} \circ \sigma =
\chi_{\sigma^{-1}(x+Z_{p_j})}.$
\item
Since $\sigma^{-1}$ takes cosets of $Z_{p_j}$
to cosets of $Z_{p_j}$, we know that $f\circ \sigma$
is a relation.
\item
Since $\sigma$ is a permutation of the
domain of $f$, the range of $f \circ \sigma$
is the same
as the range of $f$.  Thus, if $f$ is
quasi-relation (a relation with its range
a subset of $\{0, \pm 1\}$),
then so is $f \circ \sigma$.
\end{enumerate}
Hence, if $E$ is quasi-independent and $f$
is a quasi-relation on $\sigma(E)$, we would have
$f \circ \sigma$ being a quasi-relation on
$E$ and hence $f \circ \sigma =0$.   That makes $f =0$
and hence $\sigma(E)$ is quasi-independent
as well.   Likewise, if $E$ is independent and $f$ is
a relation on $\sigma(E)$, that would make
$f \circ \sigma$ be a relation on $E$ and again $f \circ
\sigma =0$.  That makes $f=0$ and $\sigma(E)$ independent.

The previous paragraph applies with
$\sigma^{-1}$ in the role of $\sigma$.
Thus, if $E$ is (qua\-si-)\-in\-de\-pend\-ent,
so is $\sigma^{-1}(E)$.   Thus, if $\sigma(E)$ is
\qiii, so is
$E=\sigma^{-1}(\sigma(E))$.

\end{proof}

\subsection{Proof of \thmref{rigidpermuts}}
\begin{proof}[Proof of
\thmref{rigidpermuts} \ \therosteritem1]
 Let $x \in Z_n$ and $1 \le j \le K$.
In the direct sum
 \begin{equation*}
Z_n = \prod _{j=1}^K Z_{p_j^{n_j}},
\]
let $x_i$ be the $i$-th coordinate of $x$.

Suppose first that $j=a$.  Consider
any $w \in x+Z_{p_a}$.  We have $w_i=x_i$ for
$i \ne a$ and $w_a \in x_a+Z_{p_a}$,
where $x_a+Z_{p_a}$ is a coset of $Z_{p_a}$
inside $Z_{p_a^{n_a}}.$
Then $\sigma(w)_i=x_i$ for $i \ne a$
and $\sigma(w)_a \in \sigma'(x_a+Z_{p_a}).$

By hypothesis, $\sigma'(x_a+Z_{p_a})=v_a+Z_{p_a}$ for
some $v_a \in Z_{p_a^{n_a}}.$   Thus we have
\[
\sigma(x+Z_{p_a}) \subset \lambda+Z_{p_a},
\]
a coset of $Z_{p_a}$ inside $Z_n$,
where $\lambda_i=x_i$ for
$i \ne a$ and $\lambda_a=v_a.$   However,
since $x+Z_{p_a}$ and $\lambda+
Z_{p_a}$ are finite sets of equal size,
and $\sigma$ is a permutation of $Z_n$,
we have
$ 
\sigma(x+Z_{p_a}) = \lambda+Z_{p_a}.
$ 

Next, suppose that $j \ne a$.  Consider
any $w \in x+Z_{p_j}$.  Then
$\sigma(w)_i = w_i$ for $i \ne a$
and $\sigma(w)_a = \sigma'(w_a).$  Let
$\lambda_i=x_i$ for $i \ne a$ and
$\lambda_a=\sigma'(w_a).$
Then $\sigma(w) \in \lambda+Z_{p_j}.$
Thus, $\sigma(x+Z_{p_j}) \subset \lambda+Z_{p_j}.$
Since these
are cosets of equal size and $\sigma$
is permutation, we have
$ 
\sigma(x+Z_{p_j}) = \lambda+Z_{p_j}.
$ 
By \lemref{permutLemma}, $\sigma$ preserves \qiind.
 \end{proof}

\begin{proof}[Proof of \thmref{rigidpermuts}
 \therosteritem2]
Let $F \subset Z_n$.  We shall prove that if $F$
is not (quasi)-in\-de\-pen\-dent, then $\sigma(F)$ is not.
This will suffice to show that $F$ is
(quasi)-independent if and only if
$\sigma(F)$ is \qiii, because
$\sigma^{-1}$ satisfies the same hypotheses
as does $\sigma$.

Suppose first that $F$ contains a full coset of $Z_2$.
Then of course, since $\sigma$
perserves the cosets of $Z_2$, $\sigma(F)$
contains the same coset of $Z_2$ and hence
is not quasi-independent (and hence not independent).

So we may assume that $F$ intersects each
coset of $Z_2$ in at most one point.
Let $f$ be a nontrivial relation that is
supported on $F$ with the appropriate range.
Let $h$ be defined as follows:
\[
h = f - \sum_{\substack{ x \in F
\\ \sigma(x) \ne x }} f(x)g_x
\]
where $g_x$ is the characteristic function
of $x+Z_2$.  By Corollary 1.2.2, each
$g_x$ is a relation.  Since $f$ is a relation, so is $h$.

We shall prove that $h$ is supported on
$\sigma(F)$ and has the appropriate range.  To do that, let $t \in Z_n.$
Suppose first that $t+Z_2$ does not intersect $F$.
Then $t$ and $\sigma(t)$ are not in $F$.
So $f(t)=0$ and $g_x(t)=0$ for $x$ in the
sum above.  So $h(t)=0$.

Next suppose that $t+Z_2$ does intersect $F$
and that $t \in F$.
If $\sigma(t) = t$, then each $g_x$ in
the sum above is $0$ at $t$ (distinct $x$ in
$F$ have disjoint cosets of $Z_2$, since $F$
meets each coset at most once).
Then $h(t)=f(t)$ with $\sigma(t)=t$.  So
$t \in \sigma(F)$, and $h(t)$ has a value from
the range of $f$.  If $\sigma(t) \ne t$,
then $h(t)=f(t)-f(t)g_t(t)=0$.

Finally, suppose that $t+Z_2$ does intersect
$F$ but $t \notin F$.  Then $f(t)=0.$
Let $w \in t+Z_2$ with $w \ne t$.  We have
$w \in F$.  If $\sigma(w)=w$, then $g_x(t)=0$
for each $g_x$ in the sum above (distinct
$x$ in $F$ have non-overlapping cosets of $Z_2$).
Hence $h(t)=0$.  If $\sigma(w) \ne w$,
then $\sigma(w)=t$ and $t \in \sigma(F)$.  Then
\[
h(t) = f(t)-f(w)g_w(t) = 0-f(w)\cdot 1=-f(w).
\]
Thus $t \in \sigma(F)$ and the value of
$h(t)$ is in $\{0, \pm 1\}$ if $f$ is a
quasi-relation, and in $\Q$ generally.

Finally, we need to observe that $h$ is nontrivial.
Since $f$ is nontrivial, there is
some $x \in F$ such that $f(x) \ne 0$.  If
$\sigma(x)=x$, then $x \in \sigma(F)$ and
(as noted above) $h(x)=f(x) \ne 0.$  Suppose
instead that $\sigma(x) \ne x$.  Let $t
=\sigma(x)$.  Then $\sigma(t)=x$ and $t \notin F$.
As noted above, $h(t)=-f(x) \ne 0.$
\end{proof}

\begin{proof}[Proof of \thmref{rigidpermuts}
\ \therosteritem3-\therosteritem6]
\therosteritem3 \ is immediate from
\thmref{primsubgroupscosets}.

\therosteritem4 \ is immediate from \therosteritem1 \ and
\thmref{primsubgroupscosets}.

\therosteritem5 \ is immediate.
Also,  \therosteritem6
 follows immediately from \therosteritem1 \ and
\therosteritem5.
 \end{proof}

 \begin{proof}[Proof of
 \thmref{rigidpermuts} \ \therosteritem7,
 $n$ square-free]
If $K=1$, there  is nothing to prove, by \thmref{rigidpermuts} \ \therosteritem1.
So assume $K\ge2.$
We will have four steps.  For all of them, we may assume
(by composing $\sigma$ with a translation)
that $\sigma(0)=0$.

\emph{(i) Reduction to $n$ is odd}.
Let $H_1=\{0\}\times Z_{p_2}
\times\cdots\times Z_{p_n}$.
Suppose that $p_1=2$.
We define a permutation $\tau$ of $Z_n$
as follows. Consider a coset
$C=\{0,1\}\times\{a\}$ of $Z_2$
with $a\in H_1$. Then $\sigma(C)=Z_2\times\{b\}$
for some $b\in H$.
Otherwise $\sigma$ would map the two element
non-quasi-independent coset
to an independent two element set.
If $\sigma(0,a)=(0,b)$,we define
$\tau(0,b)=(0,b)$ and $\tau(1,b)=(1,b)$.
Otherwise we define $\tau(0,b)=(1,b)$
and $\tau(1,b)=(0,b)$. Proceding in
this way, we define $\tau$. It's clear
that $\tau$ maps cosets of
$Z_2$ to themselves, and that $\tau\circ\sigma$
is the product of
the identity function on $Z_2$ and
some permutation on $H_1$.
It's also clear
that $\tau\circ\sigma$ preserves \qiii\ sets, by
  \thmref{rigidpermuts} \ \therosteritem2 \
  applied to $\tau$.
  We thus
may assume that $\sigma$ maps $H_1$ into
itself. We have reduced to the case that
$\sigma$
does not change the $Z_2-$coordinate of
elements of $Z_n$: that is, we have reduced
to the case of $H_1$, and
we have begun our induction (reducing $K$ by $1$)
in this case (that $p_1=2$). We have also reduced to the case of
odd $n$ (still assumed to be square-free).

For the remaining steps, we define
$
H=\zp{{1}}\times\cdots\times \zp {{K-1}}.
$

\emph{(ii) $\sigma$
maps each coset of $H$ into  another
(or the same) coset of $H$}.
We begin with consideration of $\zp1$.
Since $\sigma(0)=0$,
\propref{minsizeofqiset}
\,\therosteritem2
\ tells us that
 $\sigma(Z_{p_1})=Z_{p_1}$
because $\sigma(Z_{p_1})$ is a coset of $Z_{p_1}$
containing $0$. Hence, $\sigma$ maps $\zp1$
(which is contained in $H$) into $H$.
 Let $1\leq j<K$ and consider any coset $F$ of $\zp j$.
  Then $F$ is entirely contained in a coset of $H$.
  Suppose that $E=\sigma(F)$ meets two different cosets
  of $H$. Since $H$ has $p_K>p_j$ cosets,
  $E$ must miss at least one coset of $H$.
  Therefore, for $E$ to non \qiii \ (which it is, since
  $\sigma$ preserves the non \qiii \ sets),
  the intersection of $E$ with a coset of $H$
  must be non \qiii
  by \corref{smallqiset}.
  Since $E$ is presumed to meet at
  least 2  cosets of $H$,
  there is a proper subset $E^\prime$
  which is not \qiii. Hence
  $\sigma^{-1}(E^\prime)$ is a proper
  subset of $F$ which
  is not \qiii. But $F$ is a coset of $\zp j$, so
  every proper subset of $F$ is independent.
  This contradiction
  shows that $\sigma(F)$ is entirely
  contained in one coset
  of $H$.

Now let $m\in \zp K$ and $a,b\in H\times\{m\}$.
Then there exists a ``path'' from $a$ to $b$
via cosets of the $Z_{p_j}, 1\leq j<K$.
That is, there exists a finite set of
cosets $F_1,\dots F_s$
such that $a\in F_1$, $b\in F_s$,
$F_i\cap F_{i+1}\neq\emptyset,\, 1\leq i<s$
and each $F_i$ is a coset of one of $\zp{j}, 1\leq j<K$.
By the preceding paragraph,
there exists $\ell$ such that $\sigma(F_j)
\subset H\times\{\ell\}$
for all $j$. Therefore, $\sigma(H\times\{m\})
\subset H\times\{\ell\}$.
By cardinality, the last containment is an equality.

In the case of $H$ itself, because $\sigma(0)=0$,
we see that $\sigma(H)=H$.
By composing $\sigma $ with a permutation of
$\zp K$ that leaves $0$ fixed,
we may assume
that
\[
\sigma(H\times\{\ell\})=H\times\{\ell\}
\text{ for all }\ell.
\]

\emph{(iii) $\sigma$ maps each coset of $\zp K$ to another coset
(or the same) of $\zp K$}.
Let $F$ be a coset of $\zp K$.
Suppose that $E=\sigma(F)$ is not a coset of
$\zp K$. There are two ways this might occur. First, that
$E$ does not meet all cosets of $H$. Or, secondly, that $E$ meets all
cosets of $H$, but nevertheless, $E$ is not a coset of $\zp K$.
Now, if $E$ did not meet all cosets of $H$, then two elements of $F$
would be mapped into the same coset of $H$. Since different elements
of $F$ belong to different cosets of $H$, that's impossible, by
Step (ii). Therefore, $E$ meets every coset of $H$. Since $F$ is not
independent, $E=\sigma(F)$ is not independent.
By Remark \ref{wavyspike}, $E$ is a coset of $\zp K$.

\emph{(iv).  $\sigma$ factors.}
We have already shown that $\sigma$ maps cosets of
$\zp K$ onto cosets of $\zp K$.

It now follows that $\sigma$ is (reduced to) a
product of a permutation on $H$ with
the identity permutation of $\zp K$.
That shows that $\sigma$ is (now) determined
by what it does on $H$: that is,
we have eliminated the factor $\zp K$.
In otherwords, we have reduced $K$ by one.
By induction, we may assume that
$K=1$, which case has been taken care of.

That completes the proof of \thmref{rigidpermuts}  \therosteritem7 \ in the
 case that $n$ is square free.
\end{proof}

 \begin{proof}[Proof
 of \thmref{rigidpermuts} \ \therosteritem7:
 General case]
By applying a translation,
we see that we may assume
that $\sigma(0)=0$, as in the square free case.
We claim that
 $\sigma$
maps each coset of $Z_{\tilde n}$ onto
another (or the same) coset of $Z_{\tilde n}$.
 Indeed, suppose first that $F$ is a coset
 of some $\zp j$, $1\leq j\leq K$,
 and that $E=\sigma(F)$. If $E$ is
 not entirely contained in one
 coset of $Z_{\tilde n}$,
 then one of the intersections
 of $E$ with a coset of
 $Z_{\tilde n}$ would not be \qiii,
 by \thmref{primsubgroupscosets}. But that
 intersection
 has cardinality less than $p_j$, so
 is independent, since $\sigma$
 preserves \qiind.
 By using the path argument of
 the square-free case, we see that $\sigma$
 maps cosets of $Z_{\tilde n}$ onto
 cosets of $Z_{\tilde n}$.
 By applying \thmref{rigidpermuts} \ \therosteritem3,
 we may assume that $\sigma$ maps each coset of
 $Z_{\tilde n}$ onto itself.

  Now the proof for the case
 of square free $n$, applied to each
 coset of  $Z_{\tilde n}$,
 shows that $\sigma$ has the required form on each
 coset of $Z_{\tilde n}$, and thus
 completes the proof of
 \thmref{rigidpermuts} \  \therosteritem7.
 \end{proof}

\section{Quasi-independent set extension
  \&  $\Psi$'s
  monotonicity}\label{monotone}

\subsection{An extension theorem}
Before proving \thmref{monotonicity},
we need an extension theorem.
It shows one way to construct new quasi-independent
sets from old ones, leading to estimates for $\Psi(n)$.

\begin{theorem}\label{lambdaecupF}
Let $n>1$ be square-free, with prime factors $p_1<p_2 <
\hdots <p_K$.  Let $s \in \{1, \dots, K\}$.  If $s=K$, let $q$ be any prime such that
$q>p_K$.  If $s<K$, let $q$ be any prime such that $p_s<q$ and $q \notin
\{p_{s+1}, \hdots , p_K\}$.  Let $m=n/p_s$ and identify $Z_n$ with $Z_{p_s} \times
Z_m$; also, identify $Z_{qm}$ with $Z_q \times Z_m$.  Let $\lambda$ be the
identity embedding of $Z_n$ into $Z_{qm}$ under the
identifications
$\lambda(k,w)=(k,w)$ for $0 \le k <p_s <q$
and $w \in Z_m$.  Let $E$ be quasi-independent
in $Z_n$.  Choose $F \subset Z_{qm}$ so that
\begin{enumerate}
\item $F \cap (k+Z_m)=\emptyset$ for $k \in Z_{p_s}$.
\item For $k \ge p_s$ in $Z_q$,  $F \cap (k+Z_m)$ is  a quasi-independent set of
maximum size $\Psi(m)$.
\end{enumerate}
Then $\lambda(E) \cup F$ (or, more simply, just $E \cup F$) is quasi-independent
in $Z_{qm}$.
\end{theorem}

\begin{proof}
In $\Q[Z_{qm}]$, we use the following
basis for the kernel of $\psi$:
\begin{align}
&\text{
All characteristic functions of cosets of $Z_q$.}
\label{basiszpqr0}\\
&\text{ Any fixed set
$B_0=\{E_j\}$ of  cosets of subgroups of $Z_m$
such that}\label{basiszpqr1}
\\&\qquad\text{$\{\chi_{E_j}\}$ spans
$\text{ker\,}\psi_{|Z_m}$}\notag
\end{align}
Then $\bigcup_{1\le h<q}B_0+h$ (together with
the cosets of $Z_q$) will span the kernel of $\psi$.
That is possible by \lemref{basiszpqr}.

Suppose that $f$ were a quasirelation supported on
$R=\lambda(E)\cup F$. Then for some
rational numbers $a_t$ and $b_{h,j}$,
$ 
f=\sum_{t\in Z_{m}} a_t \chi_{t+Z_{q}} +
\sum_{\substack{j, h \in Z_q \\ 0<h<q}}
b_{h,j}\chi_{E_j+h}.
$ 
Let
$ 
g=\sum_{t\in Z_{m}} a_t \chi_{t+Z_{p_s}} +
\sum_{\substack{j,h \in Z_p \\ 0<h<p}}
b_{h,j}\chi_{E_j+h} \in \Q[Z_n].
$ 
Note that $f$ restricted to $Z_n$ is
$g$, so $g$ is a quasirelation
supported on $R\cap Z_n=E$:
that is, $g=0$.
In particular, $a_t=0 $ for all $t\in Z_m$.
Therefore
$ 
f=\sum_{\substack{j, h \in Z_q \\ 0<h<q}}
b_{h,j}\chi_{E_j+h}.
$ 
Fix $h\in Z_q, h>0$. Then
the restriction of $f$ to the coset $h+Z_m$
of $Z_m$ is a quasirelation supported
on the quasi-independent set $R\cap (h+Z_m)$. Therefore,
$
\sum_{\substack{j}}
b_{h,j}\chi_{E_j+h}=0
\text{ and so } f=0.
$
Therefore $R$ is quasi-independent.

\end{proof}

\subsection{Proof of
\thmref{monotonicity}}
\label{pfmonotonicity}
 \begin{proof}[Proof of \therosteritem1]
 The special case of the assertion for $s=K=3$
 appears as \cite[Lemma~7.1]{RGS}. The proof
 there generalizes without difficulty to the
 general case here.
 \end{proof}

 \begin{proof}[Proof of \therosteritem2]
Let $E$ be a quasi-independent set such that $\#(E)=\phi(n)+\Delta$.

We induce an imbedding of $Z_n$ into $Z_{qn/p_s}$ from the following
``lattice'' imbedding:
let $\lambda$ be the mapping of
\[
Z_{p_s} \times \prod_{j \in D \backslash\{s\}} Z_{p_j}
\text{\quad
into
\quad}
Z_q \times \prod_{j \in D \backslash\{s\}} Z_{p_j}
\]
literally (let it be the identity mapping at the level of the coordinates of
vectors).  Let
$E^*=\lambda(E) \cup F$ where $F$ is an arbitrary choice in the cosets
$w+T_m$,
for $p_s\le w <q$ and $w\in Z_q$, of
a quasi-independent set of size $\phi(m)$.
(The idea of
\cite[1.7.1]{RGS} 
suggests this approach.)

We claim that $E^*$
 has the correct cardinality, and is quasi-indepen\-dent.
For the cardinality, we compute:
\begin{align*}
\#(E^*)&=\#(E)+(q-p_s)\phi(m)
=\phi(mp_s)+\Delta+(q-p_s)\phi(m)
\\
&=(p_s-1)\phi(m)
+
(q-p_s)\phi(m)
+\Delta
=\phi(qm)+\Delta
\\
&=\phi(nq/p_s)+\Delta.
\end{align*}

To see that $E^{*}$ is quasi-independent,
apply \thmref{lambdaecupF}.
 \end{proof}

\begin{proof}[Proof of \therosteritem3]
In the proof, when expanding $E$ with $F$,
choose quasi-inde\-pen\-dent subsets of maximum size
$\Psi(n/p_s)$
in the cosets of $w+Z_{n/p_s}$ where $w \in Z_q$ and $p_s \le w <q$.
Let $E$ be a quasi-independent subset of $Z_n$ of size $\Psi(n)$.
Construct $F$
as in \thmref{lambdaecupF}, of size $(q-p_s)\Psi(m)\ge (q-p_s)\phi(m)$.
Note that $\phi(n)=(p_s-1) \phi(m)$ and $\phi(qm)=(q-1) \phi(m)$.
Because $E \cup F$ is quasi-independent, we have
\begin{align*}
\Psi(qm) &\ge \#(E \cup F)\  =\ \Psi(n)+(q-p_s)\Psi(m) \\
&\ge \Delta+
(p_s-1) \Psi(m)+(q-p_s)\Psi(m)\ =\ \Delta+(q-1) \Psi(m)
\end{align*}
\end{proof}

\subsection{Beyond $105$ and $15p$}
\label{beyondOneOhFive}
In \cite[Thm.~7.4]{RGS}  
it was shown that
\[
\Psi(105)=\phi(105)+4=52
\]
and that
\[
\Psi(15p)=\phi(15p)+4
\text{
for all primes $7\le p$.}
\]
We
can extend  those results to larger prime factors.
%

\begin{proposition}\label{psinqgemore}
 Let $n>1$  be odd with
 at least three odd  prime factors.
Then $\Psi(n)\ge \phi(n)+4.$
\end{proposition}

\begin{proof} For $n=15p$ ($p\ge7$),
$\Psi(n)=\phi(n)+4$
by \cite[Thm.~4.4]{RGA}.
Now apply \thmref{monotonicity}\ \therosteritem1
to obtain  $\Psi(n)\ge\phi(n)+4$ for
all $n$ the product of three
distinct, odd,
prime factors. For square free $n$, we  apply
\corref{psinqge}, noting that $\phi(nq)=(q-1)\phi(n)$
holds under the hypotheses of \corref{psinqge}.
Now  an application of
\thmref{Psiofn} 
\therosteritem1 \ completes the proof.
 \end{proof}

\begin{corollary}\label{psififteensmall}
\begin{enumerate}
\item  For all odd primes
$p<q<r$, $\Psi(pqr)\ge \phi(pqr)+4$.
\item  For all odd primes $p<q<r$
with $5<q$, $\Psi(pqr)\ge \phi(pqr)+5.$
\end{enumerate}
\end{corollary}
\begin{proof}
 We begin with the equality
 $\phi(15 p)+4=\Psi(15p)$
 \cite[Thm.~7.4]{RGA}.
Then
\therosteritem1 \ follows from that and
\thmref{monotonicity}.

\therosteritem2 \ follows from \propref{fiveRules}  \ and
 \thmref{monotonicity}.
\end{proof}

The time and space complexity of
the computation of $\Psi$
grows super exponentially with $n$.
However, to get lower bounds,
there is a somewhat simpler approach.

Let $n_1,\dots, n_k\ge 2$ be integers. Consider the finite abelian
group $G=Z_{n_1}\times\cdots\times Z_{n_k}$. We identify $Z_{n_j}$ with
its image $\{0\}\times\cdots\times\{0\}\times Z_{n_j}\times\{0\}\cdots$
in $G$. A {\it relation\/} on $G$ is a sum, with integer coefficients,
of characteristic functions of certain cosets of $Z_{n_j}$. Note that if $Z_{n_j} $
contains a subgroup isomophic to  $Z_{n_\ell}$,
we do {\it not\/} include characteristic function of that subgroup.
The cosets we do include are all and only those of the form
$x_1\times\cdots\times x_{j-1}\times Z_{n_j}\times x_{j+1}\cdots$.
That is, we move layers (the $Z_{n_j}$) around, but we do not do pick up a
coset inside
a layer.
A {\it quasi-relation\/} is a
relation that takes on only the values
$\pm1$ and
$0$.
A subset of $G$ is \qiii \ iff it
does not support a (quasi-)relation.

It is then not hard to see
that
\thmref{rigidpermuts}
\therosteritem1-\therosteritem4 \ extend to
this context, as does
\thmref{monotonicity}.

\begin{proposition}\label{fiveRules}
$\Psi(231)\ge \phi(231)+5$.
\end{proposition}

\begin{proof}
First, $231=3\cdot 7 \cdot 11$.
Then the proof is immediate from the extended monotonicity theorem
which tells us that
$\Psi(231)\ge \phi(231) +5$ follows from
$\Psi(3\cdot 6\cdot 9)\ge (3-1)(6-1)(9-1)+5$, and that
inequality is the assertion of the next example.
\end{proof}

An Example of $\Psi-\phi = 5$ in the $3\times6\times9$ lattice:
the following example has been verified by two independent computer programs.
It is organized by horizontal layers (9 of them),
of sizes 9, 9, 9, 9, 9, 10, 10, 10 and 10.
This totals to 85 elements, 5 more than 2(5)(8).
By layers we have
 \begin{align*}
&\{1, 2, 3, 4, 7, 8, 9, 11, 18\}
& \{1, 3, 8, 10, 12, 13, 14, 15, 17, 18\}+18
\hfill
\\
&\{3, 4, 5, 6, 8, 9, 10, 11, 13\}+36
& \{2, 3, 8, 9, 10, 11, 12, 13, 16, 17\}+54
\\
&\{1, 3, 4, 6, 9, 10, 11, 12, 14\}+72
& \{1, 2, 3, 4, 5, 7, 8, 12, 15, 17\}+90
\\
&\{1, 3, 5, 6, 7, 8, 11, 12, 16\}+108
& \{2, 3, 4, 5, 6, 7, 9, 11, 14, 1\}+126
\\
&\{1, 2, 4, 5, 7, 8, 9, 10, 18\}+144
\end{align*}

Whether $\Psi(n)-\phi(n)$ is unbounded we do not
know (even for $n$ the product of
3 primes),
much less whether $\Psi(n)/\phi(n)$ is unbounded.

When
$p=11$,
$ 
4+\phi(3\cdot 5 \cdot 11)=84>165/2.
$ 
Likewise, when $p=13$
$ 
4+\phi(3 \cdot 5 \cdot 13)=100>195/2.
$ 
So, it is possible that $Z_{165}$ and
$Z_{195}$ are the unions of two quasi-independent sets.
Whether this
is so seems to require (nearly) maximal
quasi-indepen\-dent sets that are
nearly disjoint.  See \cite[Prop.~7.5]{RGS}
for a related result.


\enddocument